\documentclass{amsart}
\usepackage{enumerate}

\def\co{\frac{\sqrt\alpha-1}{a}}
\def\Empty{\phi}
\def\x{\times}
\def\N{\mathbf N}
\def\r{\mathbf R}
\def\rn{\r^n}
\def\rm{\r^m}
\def\rnm{\rn\x\rm}
\def\rbar{[-\infty,+\infty]}

\def\set#1{\left\{#1\right\}}
\def\Ia{\set{1,...,\alpha}}

\def\ip#1{\left<#1\right>} % ip=inner product
\def\ipi#1{\ip{#1}_i}

\def\norm#1{\|#1\|}
\def\normi#1{\norm{#1}_i}
\def\normj#1{\norm{#1}_j}

\def\normmax#1{\norm{#1}_{\infty}}

\def\Lim{\lim\limits}
\def\Prod{\prod\limits}
\def\Sum{\sum\limits}
\def\Min{\min\limits}
\def\Max{\max\limits}
\def\Inf{\inf\limits}
\def\Sup{\sup\limits}
\def\d{\partial}
\def\espace{\quad}

\newtheorem{theorem}{Theorem}[section]
\newtheorem{corollary}[theorem]{Corollary}
\newtheorem{proposition}[theorem]{Proposition}
\newtheorem{lemma}[theorem]{Lemma}
\newtheorem{definition}[theorem]{Definition}

\newtheorem{remark}[theorem]{Remark}

\def\proof{\noindent{\textbf{Proof. }}}
\def\endproof{\hfill{\rule{0.5em}{0.5em}}}
\def\keywords#1\par{\\ \newline \textbf{Key words:} #1}

\begin{document}

\title{A Convergence Result for Asynchronous Algorithms and
Applications}

\author{Abdenasser BENAHMED}
\address{A. BENAHMED, Lyc\'ee Oued Eddahab Lazaret 60000 Oujda Morocco}
\email{benahmed.univ.oujda@menara.ma, nasserov@hotmail.fr}

\date{}
\maketitle
\begin{abstract}
We give in this paper a convergence result concerning parallel
asynchronous algorithm with bounded delays to solve a nonlinear
fixed point problems. This result is applied to calculate the
solution of a strongly monotone operator. Special cases of these
operators are used to solve some problems related to convex
analysis like minimization of functionals, calculus of saddle
point and variational inequality problem.\ \\
\end{abstract}

% Introduction -----------------------------------------------
\section{Introduction}\label{int}
In this paper we give a convergence result for parallel
asynchronous iterations with bounded delays. The convergence
result of these algorithms was shown by many authors. Chazan and
Miranker in \cite{cm69} treated the chaotic iterations using a
linear and contractive mapping. In 1975, Miellou\cite{mie75}
extended the works of Chazan and Miranker to the nonlinear case
using a contraction mapping and proposes a model with bounded
delays. In 1978, Baudet in \cite{bau78} generalizes the chaotic
iterations of Chazan-Miranker and Miellou and proposes a model
where the delays considered can be infinite. In a different
context, El Tarazi\cite{elt82} also established this result by a
contraction technique according to a suitable scalar norm.
Recently, Bahi\cite{bahi00} gave a convergence result concerning
parallel asynchronous algorithm to solve a linear fixed point
problems using nonexpansive linear mappings with respect to a
weighted maximum norm. Our goal is to establish a convergence
result concerning parallel asynchronous algorithm to solve a
nonlinear fixed point problems using a nonlinear and nonexpansive
mapping. All the results established in this study are mentioned
by the author in \cite{ben05}. We regard this study as a
generalization to the asynchronous case of all results stated by
Benahmed and Addou in \cite{ab05} and so, we repeat the proofs
given in \cite{ab05} by including the modifications which requires
the asynchronous case. Section \ref{pre} is devoted to a brief
description of asynchronous parallel algorithm. In section
\ref{mr} we prove the main result concerning the convergence of
the general algorithm to a fixed point of a nonlinear operator
from $\rn$ to $\rn$. This result is applied in section \ref{app}
to the operator $F=(I+cT)^{-1}$ ($c>0$) which is called the
proximal mapping associated with the maximal monotone operator
$cT$ (see Rockafellar\cite{roc76}) to calculate a solution of the
operator $T$. Special cases of these operators are also studied to
solve optimization problems and variational inequality problem.

% Section -----------------------------------------------
\section{Preliminaries}\label{pre}
$\rn$ is considered as the product space $\Prod_{i=1}^\alpha
\r^{n_i}$, where $\alpha\in\N-\{0\}$ and $n=\Sum_{i=1}^{\alpha
}{n_i}$. All vectors $x\in\rn$ considered in this study are
splitted in the form $x=(x_1,...,x_\alpha)$ where
$x_i\in\r^{n_i}$. Let $\r^{n_i}$ be equipped with the inner
product $\ipi{.,.}$ and the associated norm
$\normi{..}=\ipi{.,.}^{1/2}$. $\rn$ will be equipped with the
inner product $\ip{x,y}=\Sum_{i=1}^{\alpha }\ipi{x_i,y_i}$ where
$x,y\in\rn$ and the associated norm
$\norm{x}=\ip{x,x}^{1/2}=(\Sum_{i=1}^{\alpha}\normi{x_i}^2)^{1/2}$.
It will be equipped also with the uniform norm
$\normmax{x}=\Max_{1\leq i\leq \alpha}\normi{x_i}$.
\begin{definition}
Define $J=\set{J(p)}_{p\in\N}$ a sequence of non empty sub sets of
$\Ia$ and $S=\set{(s_1(p),...,s_\alpha(p))}_{p\in\N}$ a sequence
of $\N^\alpha$ and consider an operator
$F=(F_1,...,F_\alpha):\rn\to\rn$. The asynchronous algorithm
associated with $F$ is defined by,
\begin{equation}\label{algorithme}
\left\{
\begin{array}{l}
x^0=(x_1^0,...,x_{\alpha }^0)\in \rn \\
x_i^{p+1}=\left\{
\begin{array}{lll}
x_i^p  &  if  &  i\notin J(p)\\
F_i(x_1^{s_1(p)},...,x_{\alpha}^{s_{\alpha}(p)})  &  if  & i\in
J(p)
\end{array}
\right. \\
i=1,...,\alpha\\
p=0,1,..
\end{array}
\right.
\end{equation}
\end{definition}
It will be denoted by $(F,x^0,J,S)$. This algorithm describes the
behavior of iterative process executed asynchronously on a
parallel computer with $\alpha$ processors. At each iteration
$p+1$, the $i^{th}$ processor computes $x_i^{p+1}$ by using
(\ref{algorithme}).\\
$J(p)$ is the subset of the indexes of the components updated at the $p^{th}$ step.\\
$p-s_i(p)$ is the delay due to the $i^{th}$ processor when it
computes
the $i^{th}$ block at the $p^{th}$ iteration.\\
If we take $s_i(p)=p \ \forall i\in \Ia$, then (\ref{algorithme})
describes synchronous algorithm (without delay). During each
iteration, every processor executes a number of computations that
depend on the results of the computations of other processors in
the previous iteration. Within an iteration, each processor does
not interact with other processors, all interactions takes place
at the end of iterations.\\
If we take
$$
\left\{
\begin{array}{ll}
s_i(p)=p   &  \forall p \in \N,\forall i \in \Ia\\
J(p)=\Ia  &  \forall p \in \N
\end{array}
\right. \\
$$
then (\ref{algorithme}) describes the algorithm of Jacobi.\\
If we take
$$
\left\{
\begin{array}{ll}
s_i(p)=p   &  \forall p \in \N,\forall i \in \Ia\\
J(p)=p+1 \ (mod\ \alpha)   &  \forall p \in \N
\end{array}
\right. \\
$$
then (\ref{algorithme}) describes the algorithm of Gauss-Seidel.\\
For more details about asynchronous algorithms see \cite{cm69},
\cite{mie75}, \cite{bau78}, \cite{elt82} and \cite{bt89}.

\begin{definition}
An operator $F$ from $\rn$ to $\rn$ is said to be nonexpansive
with respect to the norm $\norm{..}$ if,
$\norm{F(x)-F(x')}\le\norm{x-x'}$ for all $x,x'\in\rn$
\end{definition}

\section{The main result}\label{mr}
We establish in this section the convergence of the general
parallel asynchronous algorithm with bounded delays to a fixed
point of a nonlinear operator $F:\rn\to\rn$.
\begin{theorem}\label{thprincipal}
Suppose\\
$(h_0)\ \exists$ a subsequence $\{p_{k}\}_{k\in \N}$ such that,
$\forall i\in\Ia,\ i\in J(p_{k})$ and $s_i(p_k)=p_k\\
(h_1)\ \exists s\in\N,$ such that,
$\forall i\in\Ia,\ \forall p\in\N,\ p-s\leq s_i(p)\leq p \\
(h_2)\ \exists u\in\rn,\ F(u)=u \\
(h_3)\ \forall x,x'\in\rn,\ \normmax{F(x)-F(x')} \le \normmax{x-x'} \\
(h_4)\ \forall x,x'\in\rn,\ \norm{F(x)-F(x')}^2 \le
\ip{F(x)-F(x'),x-x'} $\\
Then, for all $x^0\in\rn$ the sequence (\ref{algorithme}) is
convergent in $\rn$ to a fixed point $x^*$ of $F$.
\end{theorem}
\proof We follow the steps given in Addou-Benahmed \cite{ab05}
Theorem 4, with important modifications in the step \textit{(i)}.
The steps \textit{(ii)} and \textit{(iii)} are similar. We proceed
then in three steps:
\begin{enumerate}[(\it i)]
\item First, we show that the sequence $\set{\normmax{x^p-u}}_{p\in\N}$ is convergent.
For $p\in\N$, we consider the ($s+1$) iterates
$x^p,x^{p-1},...,x^{p-s}$ in the process and put
$$
z^p=\Max_{0\leq l\leq s}\normmax{x^{p-l}-u}=\Max_{p-s\leq l\leq
p}\normmax{x^l-u}
$$
Then $\forall i\in\Ia$ we have,\\
either $i\notin J(p)$ so,
$$
\begin{array}{lll}
\normi{x_i^{p+1}-u_i} & = & \normi{x_i^p-u_i}\\
  &  & \\
         & \leq & \normmax{x^p-u}\\
  &  & \\
         & \leq & \Max_{0\leq l\leq s}\normmax{x^{p-l}-u}\\
  &  & \\
         & = & z^p

\end{array}
$$

or $i\in J(p)$ so,
$$
\begin{array}{lll}
\normi{x_i^{p+1}-u_i} & = &
\normi{F_i(x_1^{s_1(p)},...,x_\alpha^{s_\alpha(p)})-F_i(u)}\\
  &  & \\
  & \leq & \normmax{F(x_1^{s_1(p)},...,x_\alpha^{s_\alpha(p)})-F(u)}\\
  &  & \\
  & \leq & \normmax{(x_1^{s_1(p)},...,x_\alpha^{s_\alpha(p)})-u}\ (\textrm{by}\ (h_3))\\
  &  & \\
  & = & \normj{x_j^{s_j(p)}-u_j}\espace(\textrm{for some j, $1\leq j\leq\alpha$})\\
  &  & \\
  & \leq & \normmax{x^{s_j(p)}-u}\\
  &  & \\
  & \leq & \Max_{p-s\leq l\leq p}\normmax{x^l-u}\espace(\textrm{use $p-s\leq s_j(p)\leq p$})\\
  &  & \\
  & = & z^p

\end{array}
$$
then
$$
\forall i\in\Ia,\ \normi{x_i^{p+1}-u_i}\leq z^p
$$
that is
$$
\normmax{x^{p+1}-u}\leq z^p
$$
therefore
$$
\begin{array}{lll}
z^{p+1} & = & \Max_{0\leq l\leq s}\normmax{x^{p+1-l}-u}\\
        & = & Max\set{\Max_{0\leq l\leq
        s-1}\normmax{x^{p-l}-u}\ ,\ \normmax{x^{p+1}-u}}\\
        & \leq & z^p
\end{array}
$$
which proves that the sequence $\set{z^p}_{p\in\N}$ is decreasing
(positive) then it's convergent. It's limit is
$$
\begin{array}{lll}
\Lim_{p\to\infty}z^p & = & \Lim_{p\to\infty}\Max_{0\leq l\leq s}
\normmax{x^{p-l}-u}\\
& &\\
& = & \Lim_{p\to\infty}\normmax{x^{p-j(p)}-u}\espace (0\leq
j(p)\leq s)\\
& &\\
& = & \Lim_{p\to\infty}\normmax{x^p-u}
\end{array}
$$
which proves that the sequence $\set{\normmax{x^p-u}}_{p\in\N}$ is
convergent and so, the sequence $\set{x^p}_{p\in\N}$ is bounded.

\item As the sequence $\set{x^{p_k}}_{k\in\N}$ is bounded
($\set{p_k}_{k\in\N}$ is defined by $(h_0)$), it contains a
subsequence noted also $\set{x^{p_k}}_{k\in\N}$ which is
convergent in $\rn$ to an $x^*$. We show that $x^*$ is a fixed
point of $F$. For this, we consider the sequence
$\{y^p=x^p-F(x^p)\}_{p\in \N}$ and prove that
$\Lim_{k\to\infty}{y^{p_k}}=0$.
$$
\begin{tabular}{lll}
$\norm{x^{p_k}-u}^2$
&  =  & $\norm{y^{p_k}+F(x^{p_k})-u}^2$\\
&  =  & $\norm{y^{p_k}}^2+\norm{F(x^{p_k})-u}^2+2\ip{F(x^{p_k})-u,y^{p_k}}$\\
\end{tabular}
$$
then
$$
\norm{y^{p_k}}^2=\norm{x^{p_k}-u}^2-\norm{F(x^{p_k})-u}^2-2\ip{F(x^{p_k})-u,y^{p_k}}
$$
however
$$
\begin{tabular}{lll}
$\ip{F(x^{p_k})-u,y^{p_k}}$
&  =  & $\ip{F(x^{p_k})-F(u),x^{p_k}-F(x^{p_k})}$\\
&  =  & $\ip{F(x^{p_k})-F(u),[x^{p_k}-F(u)]-[F(x^{p_k})-F(u)]}$\\
&  =  & $\ip{F(x^{p_k})-F(u),x^{p_k}-u}-\norm{F(x^{p_k})-F(u)}^2$\\
&  $\geq$  & $0\ (by\ (h_4))$\\
\end{tabular}
$$
so,
$$
\begin{tabular}{lll}
$\norm{y^{p_k}}^2$
&  $\le$  & $\norm{x^{p_k}-u}^2 - \norm{F(x^{p_k})-u}^2$\\
&  =       & $\norm{x^{p_k}-u}^2 - \norm{x^{p_k+1}-u}^2\ (by\ (h_0))$\\
\end{tabular}
$$
However, by (\textit{i}) the sequence
$\set{\normmax{x^p-u}}_{p\in\N}$ is convergent, then the sequence
$\set{\norm{x^p-u}}_{p\in\N}$ is also convergent with limit
$$
\begin{tabular}{lll}
$\Lim_{p\to\infty}{\norm{x^p-u}}$
&  =  & $\Lim_{k\to\infty}{\norm{x^{p_k}-u}}$\\
&  =  & $\Lim_{k\to\infty}{\norm{x^{p_k+1}-u}}$\\
&  =  & $\norm{x^*-u}$\\
\end{tabular}
$$
and so
$$
\Lim_{k\to\infty}{\norm{y^{p_k}}}=0
$$
which implies that
$$
\Lim_{k\to\infty}y^{p_k}=0
$$
and so
$$
x^*-F(x^*)=0
$$
that is $x^*$ is a fixed point of $F$.
\item We prove as in (\textit{i}) that the sequence
$\set{\normmax{ x^p-x^*}}_{p\in\N}$ is convergent, so
$$
\Lim_{p\to\infty}{\normmax{x^p-x^*}}
=\Lim_{k\to\infty}{\normmax{x^{p_k}-x^*}}=0
$$
Which proves that $x^p\to x^*$ with respect to the uniform norm
$\normmax{..}$.
\endproof
\end{enumerate}

\begin{remark}
The hypothesis $(h_0)$ means that the processors are synchronized
and all the components are infinitely updated at the same
iteration. This subsequence can be chosen by the programmer
(Bahi\cite{bahi00}).
\end{remark}
\begin{remark}
The hypothesis $(h_1)$ means that the delays dues to the
communications between processors and to the calculus are bounded,
which means that after ($s+1$) iterations, all the processors are
supposed to have update their own data (Bahi\cite{bahi00}).
\end{remark}
\begin{remark}
The hypothesis $(h_4)$ is verified by a large class of operators.
For example, the resolvent $F_\lambda=(I+\lambda T)^{-1}$ (where
$\lambda>0$) associated with a maximal monotone operator $T$ (see
Lemma \ref{lem2} below). Again, the metric projection $p_c$ of a
Hilbert space H onto a nonempty closed convex set C; that is, for
$x\in H$, $p_c(x)$ is the unique element of C which satisfies
$$\norm{x-p_c(x)}=\inf_{y\in C}\norm{x-y}$$
see for proof, Phelps\cite{phe93}, Examples 1.2.(f). In the linear
case, take for example a linear operator which is symmetric,
positive semi-definite (or simply positive) and nonexpansive, as
shown in the following proposition:
\end{remark}
\begin{proposition}
Let $A$ be a linear symmetric positive and nonexpansive operator
in $\rn$. Then $A$ verify the hypothesis $(h_4)$.
\end{proposition}
\proof Recall that an operator $A$ is said to be symmetric if for
all $x,y \in\rn,\ \ip{Ax,y}=\ip{x,Ay}$.
\begin{enumerate}[(\it i)]
\item The operator $B=I-A$ is symmetric. Indeed, $\forall x,y\in\rn$
$$
\begin{array}{lll}
\ip{Bx,y}&=&\ip{x-Ax,y}=\ip{x,y}-\ip{Ax,y}\\
         &=&\ip{x,y}-\ip{x,Ay}=\ip{x,y-Ay}\\
         &=&\ip{x,By}
\end{array}
$$
\item The operator $B$ est positive. Indeed, $\forall x\in\rn$
$$\ip{Bx,x}=\ip{x-Ax,x}=\norm{x}^2-\ip{Ax,x}\ge 0$$
since \qquad $\ip{Ax,x}\le \norm{Ax}\norm{x}\le \norm{x}^2$.
\item $A$ and $B$ are commuting operators. Indeed,
$$AB=A(I-A)=A-A^2=(I-A)A=BA$$
\item $AB$ is a symmetric operator. Indeed, $\forall x,y\in\rn$
$$
\ip{ABx,y}=\ip{Bx,Ay}=\ip{x,BAy}=\ip{x,ABy}
$$
\item $AB$ is a positive operator (see proof in \cite{gri81}, Theorem 10.7).
\item The operator $A$ verify the hypothesis $(h_4)$. Indeed, $\forall x\in\rn$
$$\ip{Ax,x}-\norm{Ax}^2=\ip{Ax,x-Ax}=\ip{Ax,Bx}=\ip{ABx,x}\ge 0$$
\end{enumerate}
\endproof

\section{Applications}\label{app}
\subsection{Solutions of maximal strongly monotone operators}\label{somm}
In this section, we apply the parallel asynchronous algorithm with
bounded delays to the proximal mapping $F=(I+cT)^{-1}$ ($c>0$)
associated with the maximal monotone operator $cT$. We say that a
multifunction $T$ from $D(T)\subseteq\rn\to\rn$ is monotone if
$$
\forall x,x'\in D(T),\ \ip{x-x',y-y'}\geq 0,\ \forall y\in Tx,\
\forall y'\in Tx'.
$$
It is said to be maximal monotone if, in addition, the graph
$$G(T)=\set{(x,y):x\in D(T)\ and\ y\in Tx}$$
is not properly contained in the graph of any other monotone
operator $T':D(T)\subseteq\rn\to\rn$. It is said to be strongly
monotone with modulus a ($a>0$) or a-strongly monotone if
$$
\forall x,x'\in D(T),\ \ip{x-x',y-y'}\geq a\norm{x-x'}^2,\ \forall
y\in Tx,\ \forall y'\in Tx'.
$$
Let $T$ be a multivalued maximal monotone operator defined from
$\rn$ to $\rn$. A fundamental problem is to determine an $x^*$ in
$\rn$ satisfying $0\in Tx^*$ which will be called a solution of
the operator $T$. The following Theorem gives a general result
concerning the solution of a maximal strongly monotone operator.

\begin{theorem}\label{thoperators}
Let $T$ be a multivalued maximal a-strongly monotone operator in
$\rn$ ($a>0$).Then
\begin{enumerate}
\item $T$ has a unique solution $x^*$.
\item Any parallel asynchronous algorithm with bounded
delays associated with the single-valued mapping $F=(I+cT)^{-1}$
where $c\geq\co$ converges in $\rn$ to the solution $x^*$ of the
problem $0\in Tx$.
\end{enumerate}
\end{theorem}
\proof We give the proof in the form of Lemmas. The two Lemmas
\ref{lem1} and \ref{lem2} were shown in \cite{ab05} by Addou and
Benahmed.
\begin{lemma}\label{lem1}(Addou-Benahmed \cite{ab05}, Theorem 4)
Let $T$ be a maximal monotone operator in $\rn$ and
$F=(I+cT)^{-1}$, $(c>0)$. Then the solutions of $T$ are exactly
the fixed points of $F$ in $\rn$.
\end{lemma}
\noindent{\textbf{Proof of Lemma \ref{lem1}.}
$$
\begin{array}{lll}
0\in Tx & \iff & x \in (I+cT)x \\
        & \iff & x = (I+cT)^{-1}x \\
        & \iff & x =Fx
\end{array}
$$
\endproof

\begin{lemma}\label{lem2}(Addou-Benahmed \cite{ab05}, Theorem 4)
Let $T$ be a maximal monotone operator and $F=(I+cT)^{-1}$
$(c>0)$. Then $F$ satisfy the hypothesis $(h_4)$.
\end{lemma}
\noindent{\textbf{Proof of Lemma \ref{lem2}.}%
Consider $x,x'$ in $\rn$ and show that
$$\norm{F(x)-F(x')}^2 \le \ip{F(x)-F(x'),x-x'}$$
Put $y=F(x)$ and $y'=F(x')$ then,
$$
\left\{
\begin{array}{l}
x\in y+cTy\\
x'\in y'+cTy'
\end{array}
\right.\\
$$
i.e.
$$
\left\{
\begin{array}{l}
x-y\in cTy\\
x'-y'\in cTy'
\end{array}
\right.\\
$$
As $c>0$, the operator $cT$ is monotone and so,
$$\ip{(x-y)-(x'-y'),y-y'} \geq 0 $$ therefore
$$\ip{x-x',y-y'} -\norm{y-y'}^2\geq 0 $$ which implies
$$\norm{F(x)-F(x')}^2 \le \ip{F(x)-F(x'),x-x'} $$
\endproof

We complete the proof of Theorem by the Lemma,
\begin{lemma}\label{lem3}
Let $T$ be a maximal a-strongly monotone operator in $\rn$ $(a>0)$
and $F=(I+cT)^{-1}$ $(c>0)$. Then
\begin{enumerate}[(a)]
\item $F$ has a unique fixed point $x^*$.
\item For $c\geq\co$, the mapping $F$ is nonexpansive with respect to the norm
$\normmax{..}$ in $\rn$.
\end{enumerate}
\end{lemma}
\noindent{\textbf{Proof of Lemma \ref{lem3}.}%
\begin{enumerate}[(a)]
\item
Consider $T'=T-aI$, $\beta=\frac{1}{1+ac}$ and $F'=(I+\beta
cT')^{-1}$. Take $x\in\rn$,
$$
\begin{array}{lll}
y=F(x) & \iff & x\in y+cTy \\
        & \iff & x\in (1+ac)y+c(T-aI)y\\
        & \iff & \beta x\in\beta(1+ac)y+\beta cT'y\\
        & \iff & \beta x\in (I+\beta cT')y\\
        & \iff & y =F'(\beta x)
\end{array}
$$
i.e., $\forall x\in\rn$, $F(x)=F'(\beta x)$.\\
As $T$ is a-strongly monotone (maximal), the operator $T'$ is
maximal monotone and then the map $F'=(I+\beta cT')^{-1}$ is
nonexpansive in $\rn$, so $\forall x,x'\in\rn$
$$\norm{F(x)-F(x')}=\norm{F'(\beta x)-F'(\beta x')}\leq\beta \norm{x-x'}$$

As $\beta=\frac{1}{1+ac}<1$, the application $F$ is contractive in
$\rn$ and then has a unique fixed point $x^*$ (Banach's fixed
point Theorem) which will be the solution of the operator $T$ by
Lemma \ref{lem1}.
\item The Euclidean norm and the uniform norm are equivalents
in $\rn$ by the relation:
$$\normmax{x}\leq\norm{x}\leq\sqrt{\alpha}\normmax{x}\ ,\ \forall
x\in\rn$$ Then, $\forall x,x'\in\rn$
$$\norm{F(x)-F(x')}\leq\beta\norm{x-x'}$$
implies
$$
\begin{array}{lll}
\normmax{F(x)-F(x')} & \leq & \beta\sqrt{\alpha}\normmax{x-x'} \\
        &  & \\
        & = & \frac{\sqrt\alpha}{1+ac}\normmax{x-x'}\\
\end{array}
$$
It is sufficient to take c such that
$$\frac{\sqrt\alpha}{1+ac}\leq 1$$
i.e.
\begin{equation}\label{c}
c\geq \co
\end{equation}
\end{enumerate}
\endproof

So, the theorem is entirely shown.

\subsection{Minimization of functional}\label{s-mf}
Let's begin by this proposition that provides a characterization
of strongly convex functions. Recall that a function $f:\rn \to
\r\cup\set{+\infty}$ is said to be strongly convexe with modulus a
($a>0$) or a-strongly convex if for all $x,x'\in\rn$ and $t\in
]0,1[$ one has
$$
f(tx+(1-t)x')\le tf(x)+(1-t)f(x')-\frac{1}{2}at(1-t)\norm{x-x'}^2
$$
The subdifferential of a proper (i.e not identically $+\infty$)
convex function $f$ on $\rn$ is the (generally multivalued)
mapping $\partial f:\rn\to\rn$ defined by
$$\partial f(x)=\set{y\in\rn\ f(x')\ge f(x) + \ip{y,x'-x},\ \forall
x'\in\rn}$$ which is a maximal monotone operator if in addition
$f$ is a lower semicontinuous function (lsc).
\begin{proposition}\label{df-fort-monotone}(Rockafellar\cite[Proposition 6]{roc76})
Let $f:\rn \to \r\cup\set{+\infty}$ be convex proper and lower
semicontinuous. Then the following conditions are equivalent:
\begin{enumerate}[(a)]
\item $f$ is a-strongly convex
\item $\d f$ is a-strongly monotone
\item whenever $y\in\d f(x)$ one has for all $x'\in\rn$:
$$f(x')\geq f(x)+\ip{y,x'-x}+\frac{1}{2}a\norm{x'-x}^2 $$
\end{enumerate}
\end{proposition}

\begin{corollary}
Let $f:\rn \to \r\cup\set{+\infty}$ be a lower semicontinuous
a-strongly convex function which is proper. Then
\begin{enumerate}
\item $f$ has a unique minimizer $x^*$.
\item Any asynchronous parallel algorithm with bounded delays associated
with the single-valued mapping $F=(I+c\partial f)^{-1}$ where
$c\geq\co$ converges to the minimizer $x^*$ of $f$ in $\rn$.
\end{enumerate}
\end{corollary}
\proof Remark that
$$
\begin{array}{lll}
0\in\d f(x) & \iff & f(x')\ge f(x)\ \forall x'\in\rn\\
              & \iff & f(x)=\Min_{x'\in\rn}{f(x')}
\end{array}
$$
so, the solutions of the operator
$\d f$ are exactly the minimizer of $f$.\\
The subdifferential $\d f$ is maximal and a-strongly monotone
(Proposition \ref{df-fort-monotone}). We apply then Theorem
\ref{thoperators} to the operator $\d f$.
\endproof

\subsection{Saddle point}\label{s-ps}
In this paragraph, we apply Theorem \ref{thoperators} to calculate
a saddle point of functional $L:\rnm \to [-\infty,+\infty]$.
Recall that a saddle point of $L$ is an element $(x^*,y^*)$ of
$\rnm$ satisfying
$$
L(x^*,y)\le L(x^*,y^*)\le L(x,y^*),\ \forall (x,y) \in \rnm
$$
which is equivalent to
$$
L(x^*,y^*)= \Inf_{x\in\rn}L(x,y^*)=\Sup_{y\in\rm}L(x^*,y)
$$
Suppose that $L(x,y)$ is proper and convex lower semicontinuous in
$x\in\rn$, concave upper semicontinuous in $y\in\rm$, then $L$ is
a proper closed saddle function in the terminology of
Rockafellar\cite{roc70a}. Let the subdifferential of $L$ at
$(x,y)\in\rnm$, $\d L(x,y)$, be defined as the set of vectors
$(z,t)\in\rnm$ satisfying
$$
\forall (x',y')\in\rnm\quad L(x,y')-\ip{y'-y,t}\le L(x,y)\le
L(x',y)-\ip{x'-x,z}
$$
then the multifunction $T_L$ defined in $\rnm$ by
$$
T_L(x,y)=\set{(z,t)\in\rnm:(z,-t)\in\d L(x,y)}
$$
is a maximal monotone operator; see
Rockafellar\cite{roc70a},\cite{roc70b}. In this case the global
saddle points of $L$ (with respect to minimizing in $x$ and
maximizing in $y$) are the elements $(x,y)$ solutions of the
problem $(0,0)\in T_L(x,y)$. That is
$$
(0,0)\in T_L(x^*,y^*) \iff (x^*,y^*)=arg\
\Min_{x\in\rn}\Max_{y\in\rm} L(x,y)
$$

\begin{definition}
The functional $L$ from $\rnm$ to $\rbar$ is said to be strongly
convex-concave with modulus a ($a>0$) or a-strongly convex-concave
if $L(x,y)$ is a-strongly convex in $x$ and a-strongly concave in
$y$.
\end{definition}

\begin{lemma}\label{lem31}
If the  functional $L$ is a-strongly convex-concave, then the
multifunction $T_L$ is an a-strongly monotone operator.
\end{lemma}
\proof Define the inner product and the norm in $\rnm$ as follows:
For $(x,y),(x',y')\in\rnm$:
$$
\left\{
\begin{array}{ll}
\ip{(x,y),(x',y')}_{\rnm}=\ip{x,x'}_{\rn}+\ip{y,y'}_{\rm}\\
\norm{(x,y)}_{\rnm}=\sqrt{\norm{x}^2_{\rn}+\norm{y}^2_{\rm}}
\end{array}
\right.
$$
we write simply as
$$
\left\{
\begin{array}{ll}
\ip{(x,y),(x',y')}=\ip{x,x'}+\ip{y,y'}\\
\norm{(x,y)}=\sqrt{\norm{x}^2+\norm{y}^2}
\end{array}
\right.
$$
Consider $(x,y),(x',y')$ in $\rnm$, $(z,t)\in T_L(x,y)$ and
$(z',t')\in T_L(x',y')$ and show that
$$
\ip{(x,y)-(x',y'),(z,t)-(z',t')}\geq a\norm{(x,y)-(x',y')}^2
$$
The function $L(x,y)$ is a-strongly convex in $x$. Proposition
\ref{df-fort-monotone} implies that the operator $\d_xL$ is
a-strongly monotone in $x$. As $z\in\d_xL(x,y)$ and
$z'\in\d_xL(x',y')$ we obtain,
$$ \ip{z-z',x-x'}\geq a\norm{x-x'}^2 $$
In the same way, $-t\in\d_yL(x,y)$, $-t'\in\d_yL(x',y')$ and
$\d_y(-L)$ is a-strongly monotone in $y$ (use proposition
\ref{df-fort-monotone} with $f(y)=-L(x,y)$) we obtain,
$$ \ip{(-t)-(-t'),y-y'}\leq -a\norm{y-y'}^2 $$
thus is
$$ \ip{t-t',y-y'}\geq a\norm{y-y'}^2 $$
therefore
$$
\ip{z-z',x-x'}+\ip{t-t',y-y'}\geq a(\norm{x-x'}^2+\norm{y-y'}^2)
$$
i.e.
$$ \ip{(z,t)-(z',t'),(x,y)-(x',y')}\geq
a\norm{(x,y)-(x',y')}^2
$$
which proves that $T_L$ is a-strongly monotone in $\rnm$.
\endproof\\

\noindent If $L$ is a  a-strongly convex-concave function and
proper closed then $T_L$ is maximal (see Rockafellar\cite{roc70b})
a-strongly monotone (Lemma \ref{lem31}). We can then apply Theorem
\ref{thoperators} to the operator $T_L$ so,
\begin{corollary}
Let $L$ be a proper closed a-strongly convex-concave function from
$\rnm$ into $[-\infty,+\infty]$. Then
\begin{enumerate}
\item $L$ has a unique saddle point $(x^*,y^*)$.
\item Any parallel asynchronous algorithm with bounded delays associated
with the single-valued mapping $F=(I+cT_L)^{-1}$ where $c\geq\co$
from $\rnm$ into $\rnm$ converges to the saddle point $(x^*,y^*)$
of $L$.
\end{enumerate}
\end{corollary}

\subsection{Variational inequality}\label{s-iv}
Let $C$ be a nonempty closed convex set in $\rn$ and $A$ a
multivalued maximal monotone operator in $\rn$ such that $D(A)=C$.
The variational inequality problem in its general form consists of
finding $x^*\in C$ satisfying
\begin{equation}\label{inegalitevariationnelle}
\exists\ y^*\in Ax^*,\ \ip{y^*,x-x^*}\geq 0,\ \forall x\in C.
\end{equation}
For $x\in\rn$, let $N_c(x)$ be the normal cone to $C$ at $x$
defined by
$$
N_c(x)=\set{y\in\rn :\ip{y,x-z}\geq 0,\ \forall z\in C}.
$$
The multifunction $T$ defined in $\rn$ by,
\begin{equation}\label{operatorn}
Tx=\left\{
\begin{array}{lll}
Ax+N_c(x) & if\ x\in C\\
\Empty    & if\ x\notin C
\end{array}
\right.
\end{equation}
is a maximal monotone operator (Rockafellar\cite{roc70c}).

\begin{lemma}\label{lem441}
If $A$ is a-strongly monotone then $T$ is an a-strongly monotone
operator.
\end{lemma}
\proof Consider $x,x'\in D(T)=C$, $y\in Tx$ and $y'\in Tx'$\\
then
$$
\left\{
\begin{array}{l}
y=y_1+y_2,\ y_1\in Ax,\ y_2\in N_c(x)\\
y'=y'_1+y'_2,\ y'_1\in Ax',\ y'_2\in N_c(x')
\end{array}
\right.
$$
however
$$
\left\{
\begin{array}{l}
y_2\in N_c(x)\Rightarrow \ip{y_2,x-z}\geq 0,\ \forall z\in C\\
y'_2\in N_c(x')\Rightarrow \ip{y'_2,x'-z}\geq 0,\ \forall z\in C
\end{array}
\right.
$$
therefore
$$
\begin{array}{lll}
\ip{y-y',x-x'} & = & \ip{y_1-y'_1,x-x'}+\ip{y_2-y'_2,x-x'}\\
              &&\\
              & = &  \underbrace{\ip{y_1-y'_1,x-x'}}_{\geq a\norm{x-x'}^2}+
              \underbrace{\ip{y_2,x-x'}}_{\geq 0}+\underbrace{\ip{y'_2,x'-x}}_{\geq 0}\\
              &&\\
              & \geq & a\norm{x-x'}^2
\end{array}
$$
\endproof
\begin{lemma}\label{lem442}
The solutions of the operator $T$ are exactly the solutions of the
variational inequality problem (\ref{inegalitevariationnelle}).
\end{lemma}
\proof
$$
\begin{array}{lll}
0\in Tx^* & \iff & 0 \in Ax^*+N_c(x^*) \\
          & \iff & \exists y^*\in Ax^* : 0\in y^*+N_c(x^*) \\
          & \iff & \exists y^*\in Ax^* : -y^*\in N_c(x^*) \\
          & \iff & \exists y^*\in Ax^* : \ip{-y^*,x^*-z}\geq 0 \ \forall z\in C\\
          & \iff & \exists y^*\in Ax^* : \ip{y^*,z-x^*}\geq 0 \ \forall z\in C\\
          & \iff & x^*\ is \ solution\ of\ (\ref{inegalitevariationnelle})
\end{array}
$$
\endproof\\
By applying the Lemma \ref{lem441}, Lemma \ref{lem442} and Theorem
\ref{thoperators}, we can write

\begin{corollary}
Let $C$ be a nonempty closed convex set in $\rn$ and $A$ a
multivalued maximal a-strongly monotone operator in $\rn$ such
that $D(A)=C$. Then
\begin{enumerate}
\item The variational inequality problem (\ref{inegalitevariationnelle})
has a unique solution $x^*$.
\item Any parallel asynchronous algorithm with bounded delays associated
with the single-valued mapping $F=(I+cT)^{-1}$ where $T$ defined
by (\ref{operatorn}) and $c\geq\co$ converges to the solution
$x^*$ of the problem (\ref{inegalitevariationnelle}).
\end{enumerate}
\end{corollary}

\end{document}